\documentclass[11pt,a4paper]{amsart} 
\usepackage[all]{xy}
\SelectTips{eu}{}
\usepackage{ifthen} 
\usepackage{amssymb}
\calclayout
\makeatletter 
\def\serieslogo@{} 
\def\@setcopyright{} 
\makeatother

\title{Rouquier's theorem on representation dimension}

\author{Henning Krause}
\address{Henning Krause\\ Institut f\"ur Mathematik\\
Universit\"at Paderborn\\ 33095 Paderborn\\ Germany.}
\email{hkrause@math.uni-paderborn.de}
\author{Dirk Kussin} 
\address{Dirk Kussin\\ Institut f\"ur Mathematik\\
Universit\"at Paderborn\\ 33095 Paderborn\\ Germany.}
\email{dirk@math.uni-paderborn.de}

\thanks{Version from September 12, 2005.}

\newtheorem{lem}{Lemma}[section]
\newtheorem{prop}[lem]{Proposition} \newtheorem{cor}[lem]{Corollary}
\newtheorem{thm}[lem]{Theorem}

\newtheorem*{Thm}{Theorem}

\theoremstyle{remark}

\newtheorem{rem}[lem]{Remark}
\newtheorem{exm}[lem]{Example}

\theoremstyle{definition}
\newtheorem{defn}[lem]{Definition}

\numberwithin{equation}{section}

\renewcommand{\mod}{\operatorname{mod}\nolimits}
\renewcommand{\leq}{\leqslant}
\renewcommand{\geq}{\geqslant}

\newcommand{\proj}{\operatorname{proj}\nolimits}

\newcommand{\add}{\operatorname{add}\nolimits}

\newcommand{\pd}{\operatorname{pd}\nolimits}

\newcommand{\gldim}{\operatorname{gl.dim}\nolimits}
\newcommand{\repdim}{\operatorname{rep.dim}\nolimits}

\newcommand{\End}{\operatorname{End}\nolimits}
\newcommand{\Hom}{\operatorname{Hom}\nolimits}

\newcommand{\Ker}{\operatorname{Ker}\nolimits}

\renewcommand{\dim}{\operatorname{dim}\nolimits}

\newcommand{\Ext}{\operatorname{Ext}\nolimits}

\newcommand{\Proj}{\operatorname{Proj}\nolimits}

\newcommand{\coh}{\operatorname{coh}\nolimits}

\newcommand{\Ab}{\mathrm{Ab}}

\newcommand{\op}{\mathrm{op}}

\newcommand{\comp}{\mathop{\raisebox{+.3ex}{\hbox{$\scriptstyle\circ$}}}}
\newcommand{\lto}{\longrightarrow}
\newcommand{\xto}{\xrightarrow}

\def\p{\phi}

\def\Ga{\Gamma}
\def\La{\Lambda}
\def\Si{\Sigma}

\def\A{{\mathcal A}}

\def\OO{{\mathcal O}}
\def\P{{\mathcal P}}

\def\S{{\mathcal S}}
\def\X{{\mathcal X}}

\def\T{{\mathcal T}}

\def\bfD{\mathbf D}
\def\bfK{\mathbf K}

\def\bfP{\mathbf P}
\def\bfR{\mathbf R}

\begin{document}

\begin{abstract}
  Based on work of Rouquier, some bounds for Aulander's representation
  dimension are discussed. More specifically, if $X$ is a reduced
  projective scheme of dimension $n$ over some field, and $T$ is a
  tilting complex of coherent $\OO_X$-modules, then the representation
  dimension of the endomorphism algebra $\End_{\OO_X}(T)$ is at least
  $n$.
\end{abstract}
\maketitle 

\section{Introduction}

Let $\La$ be an artin algebra and denote by $\mod\La$ the category of
right $\La$-modules.  Auslander defined in \cite{A} the {\em
representation dimension} of $\La$ as
$$\repdim \La=\min\{\gldim\End_\La(M)\mid\mbox{$M$ generates and
cogenerates $\mod \La$}\}.$$ It is known from work of Iyama
\cite{I} that $\repdim\La<\infty$.  Recently, Rouquier \cite{R} has
shown that there is no upper bound for the representation dimension of
an artin algebra. More precisely, given a field $k$ and an integer
$n\geq 1$, the exterior algebra $\La(k^n)$ has representation
dimension $n+1$. In this note, we use ideas from Rouquier's work and
prove the following result.

\begin{Thm}
Let $k$ be a field and $n\geq 1$ be an integer. Denote by $\La_n$
the $k$-algebra given by the following quiver with relations:
$$\xymatrix@=10pt{
{}_0\ar@/^1pc/[rr]^{x_0}_\vdots\ar@/_1pc/[rr]_{x_n}&&
{}_1\ar@/^1pc/[rr]^{x_0}_\vdots\ar@/_1pc/[rr]_{x_n}&&
{}_2&\cdots&{}_{n-1}\ar@/^1pc/[rr]^{x_0}_\vdots\ar@/_1pc/[rr]_{x_n}&&{}_n
}\quad \text{and}\quad x_ix_j=x_jx_i\;\; (0\leq i,j\leq n).$$ If $M$
is a finitely generated $\La_n$-module which generates the category of
finitely generated $\La_n$-modules, then $\gldim\End_{\La_n}(M)\geq
n$. In particular, $\repdim\La_n\geq n$.
\end{Thm}

Our aim is to give an elementary exposition of the fact that the
representation dimension is unbounded. Instead of taking the exterior
algebras, we have chosen another class of algebras in order to
simplify Rouquier's original proof. Let us stress again that most
ideas are taken from \cite{R}. However, along the way some of
Rouquier's arguments have been changed and some of the results seem to
be new.

\section{Dimensions of triangulated categories}

Let $\T$ be a triangulated category and fix subcategories
$\X,\X_1,\X_2\subseteq \T$.  Denote by $\X_1\ast\X_2$ the full
subcategory of $\T$ consisting of objects $X$ which admit an exact
triangle $X_1\to X\to X_2\to\Si X_1$ in $\T$ with $X_i\in\X_i$. Denote
by $\langle\X\rangle$ the smallest full subcategory of $\T$ which
contains $\X$ and is closed under taking coproducts, direct factors,
and all shifts. Let $\X_1\diamond\X_2=\langle\X_1\ast\X_2\rangle$.
Inductively one defines $\langle \X\rangle_0=0$ and $\langle
\X\rangle_n=\langle \X\rangle_{n-1}\diamond\langle\X\rangle$ for
$n\geq 1$. Note that the operations $\ast$ and $\diamond$ are
associative, thanks to the octahedral axiom; see \cite[Sec.~2]{BV}.

\begin{defn}[{\cite[Definition~3.1]{R}}]
Let $\T$ be a triangulated category. The {\em dimension} of $\T$ is
$$\dim\T=\min\{n\geq 0\mid\mbox{there exists $X\in\T$ such that
$\langle X\rangle_{n+1}=\T$}\}.$$
\end{defn}

\begin{lem}[{\cite[Lemma~3.3]{R}}]\label{le:dense}
Let $F\colon \S\to\T$ be an exact functor such that each object in $\T$ is
isomorphic to an object in the image of $F$. Then $\dim\S\geq\dim\T$. 
\end{lem}
\begin{proof}
If $\S=\langle X\rangle_n$, then $\T=\langle FX\rangle_n$. 
\end{proof}

\begin{lem}[{\cite[Lemma~4.11]{R}}]\label{le:comp}
Let $\T$ be a triangulated category and let
$$H_1\xto{f_1}H_2\xto{f_2}\ldots \xto{f_{n-1}}H_n\xto{f_n}H_{n+1}$$ be
a sequence of morphisms between cohomological functors $\T^\op\to\Ab$.
For each $i$, let $\X_i$ be a subcategory of $\T$ such that $f_i$
vanishes on $\X_i$ and $\X_i=\langle\X_i\rangle$. Then the composite
$f_n\comp\ldots\comp f_1$ vanishes on
$\X_1\diamond\ldots\diamond\X_n$.
\end{lem}
\begin{proof}
  Using induction, it is sufficient to prove the assertion for $n=2$. 
  Let $X_1\to X\to X_2\to\Si X_1$ be an exact triangle with
  $X_i\in\X_i$. We obtain the following commutative diagram with exact
  rows. 
$$\xymatrix{H_1X_2\ar[d]\ar[r]&H_1X\ar[d]\ar[r]&H_1X_1\ar[d]^0\\
H_2X_2\ar[d]_0\ar[r]&H_2X\ar[d]\ar[r]&H_2X_1\ar[d]\\
H_3X_2\ar[r]&H_3X\ar[r]&H_3X_1}$$ A simple diagram chase shows that
$f_2\comp f_1$ vanishes on $X$.  Now observe that $\X_1\diamond\X_2$
consists of direct factors of objects in
$\X_1\ast\X_2$. 
\end{proof}

Let $\A$ be an abelian category. We denote by $\bfD^b(\A)$ the derived
category of bounded complexes in $\A$ and identify $\A$ with the full
subcategory consisting of complexes concentrated in degree zero.  The
category $\bfD^b(\A)$ is obtained from the homotopy category
$\bfK^b(\A)$ by formally inverting all quasi-isomorphisms.  Note that
$$\Ext^n_\A(X,Y)\cong\Hom_{\bfD^b(\A)}(X,\Si^nY)$$for all $X,Y\in\A$
and $n\geq 0$. 

\begin{lem}\label{le:pd}
Let $\A$ be an abelian category and $X\in\A$ satisfying
$\Ext_\A^n(X,-)\neq 0$. Then we have
$X\not\in\langle\P\rangle_n\subseteq\bfD^b(\A)$ for the subcategory
$\P\subseteq\A$ of projective objects. 
\end{lem}
\begin{proof}
Let $X=X_0$ and fix an extension $$\xi\colon\;0\to X_n\to E_n\to
E_{n-1}\to\ldots\to E_1\to X_0\to 0$$ in $\A$. We view $\xi$ in
$\Ext_\A^n(X_0,X_n)$ as the composite of extensions
$$\xi_i\colon\;0\to X_i\to E_i\to X_{i-1}\to 0.$$ Each $\xi_i$ induces
a connecting morphism
$$f_i\colon\Hom_{\bfD^b(\A)}(-,\Si^{i-1}X_{i-1})\lto
\Hom_{\bfD^b(\A)}(-,\Si^{i}X_{i})$$ vanishing on $\langle\P\rangle$. 
The composite $f_n\comp\ldots\comp f_1$ sends the identity morphism of
${X_0}$ to $\xi=\xi_n\comp\ldots \comp\xi_1$. Thus $\xi\neq 0$ implies
$X_0\not\in\langle\P\rangle_n$, by Lemma~\ref{le:comp}. 
\end{proof}

\begin{lem}[{\cite[Theorem~8.3]{C}}]\label{le:C}
Let $\A$ be an abelian category and suppose that $\A$ has enough
projective objects. Let $X$ be a bounded complex in $\A$ such that
$B_iX$ and $H_iX$ have projective dimension at most $n$ for all
$i$. Then $X\in\langle\P\rangle_{n+1}\subseteq\bfD^b(\A)$ for the
subcategory $\P\subseteq\A$ of projective objects.
\end{lem}
\begin{proof}
For each $i$, choose epimorphisms $P^{B_iX}\to B_iX$ and
$P^{H_iX}\to H_iX$ such that $P^{B_iX}$ and $P^{H_iX}$ are projective. 
We obtain commutative diagrams
$$\xymatrix{
0\ar[r]&P^{B_iX}\ar@{>>}[d]\ar[r]&P^{Z_iX}\ar@{>>}[d]\ar[r]&P^{H_iX}\ar@{>>}[d]\ar[r]&0\\
0\ar[r]&{B_iX}\ar[r]&{Z_iX}\ar[r]&{H_iX}\ar[r]&0}$$
with exact rows. Similarly, we obtain  commutative diagrams
$$\xymatrix{
0\ar[r]&P^{Z_iX}\ar@{>>}[d]\ar[r]&P^{X_i}\ar@{>>}[d]\ar[r]&P^{B_{i-1}X}\ar@{>>}[d]\ar[r]&0\\
0\ar[r]&{Z_iX}\ar[r]&{X_i}\ar[r]&{B_{i-1}X}\ar[r]&0}$$ with exact
rows.  Defining $P_i=P^{X_i}$ and taking the composite
$$P_{i+1}=P^{X_{i+1}}\to P^{B_iX}\to P^{Z_iX}\to P^{X_i}=P_i$$ as
differential, we obtain a complex $P$. Note that $B_iP=P^{B_iX}$,
$Z_iP=P^{Z_iX}$, $H_iP=P^{H_iX}$ for all $i$, and therefore
$P\in\langle\P\rangle$.  Now let $X^0=X$ and $P^0=P$. We have an
epimorphism $P^0\to X^0$ and denote by $X^1$ the shift of the
degreewise kernel. This yields an exact triangle $P^0\to X^0\to
X^1\to\Si P^0$ in $\bfD^b(\A)$. Note that $B_iX^1$ and $H_iX^1$ have
projective dimension at most $n-1$ for all $i$.  We inductively
continue this construction, and by our assumption on $B_iX$ and $H_iX$, we
have $X^n\in\langle\P\rangle$. Thus $X^{n-1}\in\langle\P\rangle_{2}$,
and inductively, $X\in\langle\P\rangle_{n+1}$.
\end{proof}

Recall that a ring $\La$ is {\em right coherent} if the category
$\mod\La$ of finitely presented right $\La$-modules is abelian. 

\begin{prop}\label{pr:coh}
Let $\La$ be a right coherent ring. Then 
$$\dim\bfD^b(\mod \La)\leq \sup\{\pd X\mid X\in\mod\La\}\leq \gldim\La.$$
\end{prop}
\begin{proof} 
We apply Lemma~\ref{le:C}. Take the abelian category $\A=\mod\La$ and
observe that $\langle\P\rangle_n=\langle\La\rangle_n$ for all $n\geq 0$. 
\end{proof}

\begin{rem} 
The result of the preceding proposition answers a question in
\cite[Rem.~7.27]{R}. 
\end{rem}

\section{Representation generators}

Given an object $M$ of an additive category, we denote by $\add M$ the
smallest full subcategory closed under finite coproducts and direct
factors. 

\begin{defn}
  Let $\A$ be an abelian category and fix an object $M\in\A$. An {\em
    $M$-resolution} of an object $X\in\A$ is an exact sequence
$$\ldots\to M_2\to M_1\to M_0\to X\to 0$$ in $\A$ such that
$M_i\in\add M$ for all $i$ and the induced sequence
$$\ldots\to\Hom_\A(M,M_2)\to\Hom_\A(M,M_1)\to \Hom_\A(M,M_0)\to
\Hom_\A(M,X)\to 0$$ of abelian groups is exact.  The resolution has
{\em finite length} if $M_i=0$ for $i\gg 0$.  We call $M$ a {\em
representation generator} if every object $X\in\A$ admits an
$M$-resolution of finite length. 
\end{defn}

\begin{rem}
  Note that we do not assume any common bound for the length of an
  $M$-resolution in our definition of a representation generator $M$. 
  However, the finite bound $\gldim\End_\A(M)$ is automatic in some
  interesting cases; see Lemma~\ref{le:noethrepgen}. 
\end{rem}

\begin{lem}\label{le:adj}
Let $\A$ be an abelian category. An object $M\in\A$ is a
representation generator if and only if the inclusion $\bfK^b(\add
M)\to\bfK^b(\A)$ admits a right adjoint $F\colon
\bfK^b(\A)\to\bfK^b(\add M)$ such that the adjunction morphism $FX\to
X$ is a quasi-isomorphism for all $X\in\bfK^b(\A)$. 
\end{lem}
\begin{proof}[Proof {\rm (cf.\ \cite[Proposition~8.3]{R})}] 
Let $M$ be a representation
generator. For each $X\in\bfK^b(\A)$, we need to construct
an {\em approximation} $M^X\to X$ such that 
\begin{enumerate}
\item $M^X\in\bfK^b(\add M)$, 
\item the induced morphism
$\Hom_{\bfK^b(\add M)}(Y,M^X)\to\Hom_{\bfK^b(\A)}(Y,X)$ is bijective for
all $Y\in\bfK^b(\add M)$, and 
\item $M^X\to X$ is a quasi-isomorphism. 
\end{enumerate}
Then we define $F\colon \bfK^b(\A)\to\bfK^b(\add
M)$ by sending $X$ to $M^X$.  We construct $M^X\to X$ by induction on
the width of $X$. Suppose first that $X$ is concentrated in degree
zero. Take a finite length $M$-resolution
$$0\to M_n\to\ldots\to M_1\to M_0\to X\to 0$$
which exists by
assumption, and define $M^X_i=M_i$ for $0\leq i\leq n$ and $M^X_i=0$
otherwise. The morphism $M_0\to X$ induces a quasi-isomorphism
$\p\colon M^X\to X$. Moreover, $\Hom_{\A}(M,\p)$ is a
quasi-isomorphism and therefore $\Hom_{\bfK^b(\A)}(M,\p)$ is
bijective. Using induction on the width of a complex $Y$ in
$\bfK^b(\add M)$, one sees that $\Hom_{\bfK^b(\A)}(Y,\p)$ is
bijective. Now suppose that $X$ fits into an exact triangle $X'\to
X''\to X\to\Si X'$ where approximations $\p'\colon M^{X'}\to X'$ and
$\p''\colon M^{X''}\to X''$ have been constructed. Using the bijection
$$\Hom_{\bfK^b(\add M)}(M^{X'},M^{X''})\to\Hom_{\bfK^b(\A)}(M^{X'},X''),$$
we obtain a
morphism $M^{X'}\to M^{X''}$ and complete it to an exact triangle
$$M^{X'}\to M^{X''}\to M^X\to \Si M^{X'}.$$ Moreover, we can complete
$\p'$ and $\p''$ to a morphism of triangles and obtain an
approximation $\p\colon M^X\to X$. 

Suppose now that the inclusion $\bfK^b(\add M)\to\bfK^b(\A)$ admits a
right adjoint $F$ such that the adjunction morphism $FX\to X$ is a
quasi-isomorphism for all $X$ in $\bfK^b(\A)$. Given an object
$X\in\A$, we view $X$ as a complex concentrated in degree zero.  The
properties of the morphism $FX\to X$ imply that we obtain an
$M$-resolution
$$0\to M_n\to\ldots\to M_1\to M_0\to X\to 0$$ by taking 
$$M_i=\begin{cases}(FX)_i&i>0,\\
Z_0(FX)&i=0,\\
0&i<0. 
\end{cases}$$ 
It is easily checked that $M_0\in\add M$, and $FX\to X$ induces the
morphism $M_0\to X$.
\end{proof}

\begin{prop}\label{pr:repgen}
Let $\A$ be an abelian category with a representation generator $M$. 
Then $\Ga=\End_\A(M)$ is right coherent and every finitely presented
$\Ga$-module has finite projective dimension. Moreover, we have
$$\dim\bfD^b(\A)\leq \dim\bfD^b(\mod \Ga)\leq \sup\{\pd X\mid
X\in\mod\Ga\}\leq \gldim\Ga.$$
\end{prop}
\begin{proof}
  We use that $\Hom_\A(M,-)$ induces an equivalence $\add M\to
  \proj\Ga$ onto the category of finitely generated projective
  $\Ga$-modules.  Recall that $\Ga$ is right coherent if the kernel of
  any morphism $\p$ in $\proj\Ga$ is finitely presented. We have
  $\p=\Hom_\A(M,\psi)$ for some morphism $\psi$ in $\add M$. Then
  $X=\Ker\psi$ admits a finite length $M$-resolution $$0\to
  M_n\to\ldots\to M_1\to M_0\to X\to 0$$ by our assumption on $M$. The
  morphism $M_1\to M_0$ induces a projective presentation
$$\Hom_\A(M,M_1)\to\Hom_\A(M,M_0)\to\Hom_\A(M,X)\to 0$$ of
$\Hom_\A(M,X)=\Ker\p$ in $\mod\Ga$. Thus $\Ga$ is right coherent. 
Now let $Z$ be in $\mod\Ga$ with projective presentation
$$P_1\xto{\p} P_0\to Z\to 0.$$
Again, we have $\p=\Hom_\A(M,\psi)$ for
some morphism $\psi$ in $\add M$, and $X=\Ker\psi$ admits an
$M$-resolution $$0\to M_n\to\ldots\to M_1\to M_0\to X\to 0.$$
Applying
$\Hom_\A(M,-)$, we obtain a projective resolution of length $n$ of
$\Hom_\A(M,X)=\Ker\p$. Thus $\pd Z\leq n+2$. 

For the bounds on $\dim\bfD^b(\A)$, observe that we have exact
equivalences
$$\bfK^b(\add M)\xto{\sim}\bfK^b(\proj\Ga)\xto{\sim}\bfD^b(\mod\Ga).$$
The functor $\bfK^b(\add M)\to\bfD^b(\A)$ is essentially surjective by
Lemma~\ref{le:adj}, and we obtain
$$\dim\bfD^b(\A)\leq \dim\bfD^b(\mod \Ga)$$ from
Lemma~\ref{le:dense}. For the rest, apply Proposition~\ref{pr:coh}. 
\end{proof}

For the module category of an artin algebra, we have the following
characterization of a representation generator. 

\begin{lem}\label{le:noethrepgen}
  Let $\La$ be an artin algebra. Then $M\in\mod\La$ is a
  representation generator of $\mod\La$ if and only if $M$ generates
  $\mod\La$ and $\End_\La(M)$ has finite global dimension. 
\end{lem}
\begin{proof}
  Suppose first that $M$ is a representation generator of $\mod\La$
  with $\Ga=\End_\La(M)$. For each $X\in\mod\La$, we have an
  epimorphism $M_0\to X$ with $M_0$ in $\add M$. Thus $M$ generates
  $\mod\La$.  Every finitely presented $\Ga$-module has finite
  projective dimension, by Proposition~\ref{pr:repgen}.  It follows
  that $\Ga$ has finite global dimension, because the global dimension
  equals the projective dimension of $\Ga/\mathfrak r$, where
  $\mathfrak r$ denotes the Jacobson radical.
  
  Now suppose that $M$ generates $\mod\La$ and $\Ga=\End_\La(M)$ has
  finite global dimension. Every finitely presented $\La$-module $X$
  admits an epimorphism $\p\colon M^X\to X$ with $M^X$ in $\add M$
  such that $\Hom_\La(M,\p)$ is an epimorphism. This is clear since
  $\Hom_\La(M,X)$ is finitely generated over $\Ga$. Let $M_0=M^X$ and
  $X_1=\Ker\p$. Inductively, we define $M_i=M^{X_i}$ and obtain
  morphisms $M_i\to X_i\to M_{i-1}$ which induce a projective
  resolution
$$\ldots\to\Hom_\La(M,M_1)\to\Hom_\La(M,M_0)\to\Hom_\La(M,X)\to 0$$ in
$\mod\Ga$. We have $X_i\in\add M$ for $i\geq \gldim\Ga$ and therefore
the construction terminates, giving an $M$-resolution of finite length
$$0\to M_n\to\ldots\to M_1\to M_0\to X\to 0$$ as required. 
\end{proof}

\begin{cor}[{\cite[Proposition~8.3]{R}}]\label{co:gen}
  Let $\La$ be an artin algebra and $M\in\mod\La$.  If $M$ generates
  $\mod\La$, then
$$\dim\bfD^ b(\mod\La)\leq\gldim\End_\La(M).$$
\end{cor}
\begin{proof}
  Combine Lemma~\ref{le:noethrepgen} and Proposition~\ref{pr:repgen}.
\end{proof}

\begin{exm}[{\cite[Theorem~10.2]{A1}}]
Let $\La$ be an artin algebra with Jacobson radical $\mathfrak r$ and
$\mathfrak r^n=0$. Then $M=\coprod_{i=1}^n\La/\mathfrak r^i$ is a
representation generator of $\mod\La$ with $\gldim\End_\La(M)\leq
n$. 
\end{exm}

\section{Some geometry}

We compute explicit bounds for the dimension of certain derived
categories, using some elementary facts from algebraic geometry. Given
a scheme $X$, the category of coherent $\OO_X$-modules is denoted by
$\coh X$.

\begin{lem}\label{le:closed1}
Let $X$ be a reduced projective scheme over a field, and let $F\in\coh
X$. Then every irreducible component of $X$ contains a closed point
$x\in X$ such that $F_x$ is a free $\mathcal{O}_x$-module.
\end{lem}
\begin{proof}
Assume first that $X$ is irreducible. The local ring $\OO_\xi$ at the
  generic point $\xi\in X$ is a field, since $X$ is reduced. Thus
  $F_{\xi}$ is a free $\mathcal{O}_{\xi}$-module, and there is a
  neighbourhood $U$ such that $F|_{U}$ is free \cite[Ex.~II.5.7]{H}.
  We find a closed $x\in U$ since the closed points of $X$ are dense
  \cite[Ex.~II.3.14]{H}. 

Now assume that $X$ is arbitrary and fix an irreducible component
$Y\subseteq X$. We pass to $Y$ and the same argument as before works
if we choose the point $x\in Y$ such that it is not contained in any
other irreducible component.
\end{proof}

\begin{lem}\label{le:closed2}
  Let $X$ be a reduced projective scheme over a field, and
  $M\in\mathbf{D}^b (\coh X)$. Then every irreducible component of $X$
  contains a closed point $x\in X$ such that
  $$M_x\in\langle{\mathcal{O}_x}\rangle\subseteq\bfD^b(\mod\mathcal{O}_x).$$
\end{lem}
\begin{proof}
We apply Lemma~\ref{le:closed1} and find a closed point $x$ such that $B_i
  M_x$ and $H_i M_x$ are free $\OO_x$-modules for all $i$. Thus
  $M_x\in\langle{\mathcal{O}_x}\rangle$ by Lemma~\ref{le:C}. 
\end{proof}

\begin{lem}[{\cite[Lemma~7.14]{R}}]\label{le:local}
Let $\La$ be a commutative local noetherian ring with maximal ideal
$\mathfrak m$. If $\La$ has Krull dimension $n$, then $\La/\mathfrak
m\not\in\langle\La\rangle_n$. 
\end{lem}
\begin{proof} 
The projective dimension of $\La/\mathfrak m$ is at least $n$. Now apply
Lemma~\ref{le:pd}. 
\end{proof}

\begin{prop}[{\cite[Proposition~7.17]{R}}]\label{pr:Pn}
Let $X$ be a reduced projective scheme over a field. 
Then $\dim\bfD^b(\coh X)\geq \dim X$. 
\end{prop}
\begin{proof}
  We may identify $X=\Proj S$ for some graded ring $S$
  \cite[Cor.~II.5.16]{H}.  Let $M\in\mathbf{D}^b (\coh X)$ and $n\geq
  0$ such that $\mathbf{D}^b (\coh X)= \langle M\rangle_{n+1}$.
  Choose an irreducible component $Y\subseteq X$ of maximal dimension.
  Using Lemma~\ref{le:closed2}, there exists a closed $x\in Y$ such
  that $M_x \in\langle\mathcal{O}_x\rangle$. Let $I(x)\subseteq S$ be
  the homogeneous prime ideal corresponding to $x$ and $F$ be the
  sheaf associated to $S/I(x)$. Then $F\in\langle M\rangle_{n+1}$
  implies for the residue field
  $$k(x)=F_x\in\langle
  M_x\rangle_{n+1}\subseteq\langle\mathcal{O}_x\rangle_{n+1}.$$ Thus
  $n\geq \dim\OO_x=\dim Y=\dim X$, by Lemma~\ref{le:local}. 
\end{proof}

\section{Representation dimensions}

We are now in a position to give the proof of the main result.  

\begin{proof}[Proof of the Main Theorem.]
  Let $k$ be a field and $n\geq 1$ an integer. We consider the
  projective scheme $X=\bfP_k^n$ and fix the sheaf
  $T=\coprod_{i=0}^n\OO_X(i)$. We have
  $\End_{\OO_X}(T)\cong\La_n$ and an equivalence
$$\bfR\Hom_{\OO_X}(T,-)\colon \bfD^b(\coh
X)\xto{\sim}\bfD^b(\mod\La_n)$$ by work of Beilinson \cite{B}, because
$T$ is a tilting sheaf. Now let $M$ be a finitely generated
$\La_n$-module which generates $\mod\La_n$.  Applying
Proposition~\ref{pr:Pn} and Corollary~\ref{co:gen}, we obtain
$$n=\dim X\leq\dim\bfD^b(\coh X)=\dim\bfD^b(\mod\La_n)\leq
\gldim\End_\La(M).$$
From the definition of the representation
dimension, it follows that $\repdim\La_n\geq n$.
\end{proof}

\begin{rem} 
We have actually $\dim\bfD^b(\coh X)=n$ for $X=\bfP^n_k$ since
$$\dim\bfD^b(\coh X)=\dim\bfD^b(\mod\La_n)\leq \gldim\La_n=n$$ by
Proposition~\ref{pr:coh}. 
\end{rem}

The proof of the main theorem provides the following method for
constructing algebras with representation dimension bounded from
below. 

\begin{thm}
Let $X$ be a reduced projective scheme over a field, and let $T$ be a
tilting complex in $\bfD^b(\coh X)$ with $\La=\End_{\bfD^b(\coh
X)}(T)$.  Given a generator $M$ of $\mod\La$, we have
$\gldim\End_\La(M)\geq \dim X$. In particular, $\repdim\La\geq \dim
X$. 
\end{thm}

The actual computation of the representation dimension of
$\La=\End_{\bfD^b(\coh X)}(T)$ seems to be difficult. However, there
is the following bound from above. The argument for the first part is
due to Lenzing. We fix a field $k$ and let  $D=\Hom_k(-,k)$.

\begin{prop}
Let $X$ be a smooth projective scheme of dimension $n\geq 1$ over a field,
and let $T\in\coh X$ be a tilting sheaf with $\La=\End_{\OO_X}(T)$. 
\begin{enumerate}
\item Let $Q$ be an injective and $P$ be a projective
$\La$-module. Then $\Hom_\La(Q,P)=0$. 
\item We have $\repdim\La\leq 2\gldim\La + 1<\infty$. 
\end{enumerate}
\end{prop}
\begin{proof}
(1) Let $X$ be projective over the field $k$ and denote by
$\omega=\omega_{X/k}$ the dualizing sheaf. Then Serre duality gives
$$D\Hom_{\bfD^b(\coh X)}(E,F)\cong \Hom_{\bfD^b(\coh
X)}(F,E\otimes_{\OO_X}\omega[n])\quad\text{for}\quad E,F\in\bfD^b(\coh X)\,.$$
The equivalence
$$\bfR\Hom_{\OO_X}(T,-)\colon\bfD^b(\coh X)\xto{\sim}\bfD^b(\mod\La)$$
identifies the full subcategory
$$\A=\{F\in\bfD^b(\coh X)\mid\Hom_{\bfD^b(\coh X)}(T,F[i])=0\text{ for
}i\neq 0\}$$ with $\mod\La$.  Clearly, $T$ is a projective generator
of $\A$ and we claim that $U=T\otimes_{\OO_X}\omega[n]$ is an injective
cogenerator of $\A$. In deed, we have $U\in\A$ since
$$\Hom_{\bfD^b(\coh X)}(T,U[i])\cong D\Hom_{\bfD^b(\coh
X)}(T,T[-i])=0$$ for $i\neq 0$. Moreover, $U$ is an
injective cogenerator of $\A$ since Serre duality implies 
$$\Hom_\A(F,U)\cong D\Hom_\A(T,F)\quad\text{for}\quad F\in\A$$
and $T$ is a projective generator. Finally, we have
$$\Hom_\A(U,T)=\Hom_{\bfD^b(\coh
X)}(T\otimes_{\OO_X}\omega[n],T)\cong\Ext_X^{-n}(T\otimes_{\OO_X}\omega,T)=0\, .$$

(2) Let $M=\La\amalg D\La$. Clearly, $M$ generates and cogenerates
    $\mod\La$. Using that $\Hom_\La(D\La,\La)=0$, we obtain
$$\repdim\La\leq\gldim\End_\La(M)\leq 2\gldim\La + 1,$$ for instance
by \cite[Prop.~7.5.1]{MR}. Finally, observe that
$\gldim\La<\infty$, since $\bfD^b(\mod\La)$ admits a Serre functor. 
This follows from Lemma~\ref{le:gldim} below. 
\end{proof}

\begin{lem}\label{le:gldim}
Let $\La$ be a finite dimensional algebra over a field $k$. Suppose
there is an exact functor $F\colon\bfD^b(\mod\La)\to\bfD^b(\mod\La)$
such that 
$$D\Hom_{\bfD^b(\mod\La)}(X,Y)\cong
\Hom_{\bfD^b(\mod\La)}(Y,FX)\quad\text{for}\quad
X,Y\in\bfD^b(\mod\La).$$ Then $\gldim\La<\infty$. 
\end{lem}
\begin{proof}
Denote by $S_1,S_2,\ldots,S_r$ the simple $\La$-modules. We need to check
that for each pair $i,j$, we have $\Ext_\La^l(S_i,S_j)=0$ for $l\gg 0$. 
We compute 
$$D\Ext_\La^l(S_i,S_j)\cong D\Hom_{\bfD^b(\mod\La)}(S_i,S_j[l])\cong
\Hom_{\bfD^b(\mod\La)}(S_j,FS_i[-l]).$$ Now use that for any pair
$X,Y$ of bounded complexes, we have
$$\Hom_{\bfD^b(\mod\La)}(X,Y[-l])=0\quad\text{for}\quad l\gg 0.$$
\end{proof}

\subsection*{Acknowledgement} 
Amnon Neeman pointed out the relevance of Christensen's work
\cite{C}. We are grateful to him and to Helmut Lenzing for discussions
on the geometric aspects of the representation dimension.

\end{document}